\documentstyle[11pt]{article}
\def\R{\ifmmode{\rm I\mkern-3.1mu
R\mkern1mu}\else{\rm I\kern-.18em  R\hskip1pt\ 
}\fi\relax}  
\font\bigreek=Symbol at 16pt
\font\smgreek=Symbol at 10pt
 
\def\b{\beta}

\def\l{\lambda}
\def\L{\Lambda}
\def\n{\nu}

\def\ph{\phi}
\def\Ph{\Phi}

\def\s{\sigma}
\def\sou{\overline}
\def\so{\underline} 

\def\f{\rightarrow}
\def\q{\forall}
\def\e{\exists}

\def\p{\succ}

\def\R{\ifmmode{\rm I\mkern-3.1mu
R\mkern1mu}\else{\rm I\kern-.18em 
R\hskip1pt\ }\fi\relax} 

\def\Z{\ifmmode{ Z\mkern-4.6mu
Z\mkern2mu}\else{ Z\kern-.28em 
Z\hskip1pt\ }\fi\relax} 

\def\Q{\ifmmode{\rm Q\mkern-10mu
l\mkern4.5mu}\else{\rm Q\kern-.57em
l\hskip3pt\ }\fi\relax} 

\def\N{\ifmmode{\rm I\mkern-3.1mu
N\mkern0.5mu}\else{\rm I\kern-.16em
N\hskip0.5pt\ }\fi\relax} 

\def\C{\ifmmode{\rm C\mkern-8.8mu
l\mkern4mu}\else{\rm C\kern-.48em
l\hskip2.6pt\ }\fi\relax} 

\def\mats{\ifmmode{ {\hbox{\bigreek s}} }\else{ 
{\bigreek s} }\fi\relax}
\def\matsin{\ifmmode{ {\hbox{\smgreek s}} }\else{ 
{\smgreek s} }\fi\relax}
\def\matt{\ifmmode{ {\hbox{\bigreek t}} }\else{ 
{\bigreek t} }\fi\relax}
\def\mattin{\ifmmode{ {\hbox{\smgreek t}} }\else{ 
{\smgreek t} }\fi\relax}

\begin{document}
\begin{center}

{\Large \bf A CONJECTURE ON NUMERAL SYSTEMS}\\ [1cm]

{\bf  Karim NOUR} \\ 

LAMA - \'Equipe de Logique \\ 
Universit\'e de Savoie \\
73376 Le Bourget du Lac \\
FRANCE \\
{\it E-mail}: nour@univ-savoie.fr \\[1,5cm]

\end{center}

\begin{abstract} A numeral system is  an infinite sequence of different closed normal
$\l$-terms intended to code the integers in $\l$-calculus.  H. Barendregt has shown that if we
can represent, for a numeral system, the functions : Successor,  Predecessor, and Zero Test,
then all total recursive functions can be represented. In this paper we prove the independancy
of these particular three functions. We give at the end a conjecture on the number of unary
functions necessary to represent all total recursive functions. \\  \end{abstract}

\section{Introduction}

A numeral system is an infinite sequence of different closed $\b \eta$-normal $\l$-terms $\bf
d$ \rm = $d_0 , d_1 ,..., d_n , ...$ intended to code the integers in $\l$-calculus. \\

For each numeral system $\bf d$, we can represent total numeric functions as follows :\\

A total numeric function $\phi : \N^p \f \N$ is $\l$-definable with respect to {\bf d} iff
\begin{center} $\e$ $F_{\phi}$ $\q$ $n_1,...,n_p \in \N$ $(F_{\phi} ~ d_{n_1}... d_{n_p})
\simeq\sb{\b} d_{\phi(n_1,...,n_p)}$  \end{center}

One of the differences between our numeral system definition and the H. Barendregt's definition
given in [1] is the fact that the $\l$-terms $d_i$ are normal and different. The last
conditions allow with some fixed reduction strategies (for example the left reduction
strategy) to find the exact value of a function computed on arguments.\\

H. Barendregt has shown that if we can represent, for a numeral system, the functions :
Successor,  Predecessor, and Zero Test, then all total recursive functions can be
represented.\\

We prove in this paper that this three particular functions are independent. We think it is,
at least, necessary to have three unary functions to represent all total recursive functions.
\\

This paper is organized as followsÊ:
\begin{itemize}
\item The section 2 is devoted to preliminaries.
\item In section 3, we define the numeral systems, and we present the result of H. Barendregt.
\item In section 4, we prove the independancy of the functions: Successor,  Predecessor, and
Zero Test. We give at the end a conjecture on the number of unary functions
necessary to represent all total recursive functions.  

\end{itemize}

\section{Notations and definitions}

The notations are standard (see [1] and [2]).

\begin{itemize}
\item We denote by $I$ (for Identity) the $\l$-term $\l x x$, $T$ (for True) the $\l$-term
$\l x \l y x$ and by $F$ (for False) the $\l$-term $\l x \l y y$.
\item The pair $<M,N>$ denotes the $\l$-term $\l x (x ~ M ~ N)$.
\item The {\it $\b$-equivalence} relation is denoted by $M \simeq\sb{\b} N$. 
 \item The notation $\s(M)$ represents the result of the simultaneous
substitution $\s$ to the free variables of $M$ after a suitable renaming of the bound
variables of $M$.   
 \item A {\it $\b \eta$-normal} $\l$-term is a $\l$-term which does not contain neither a
$\b$-redex [i.e. a $\l$-term of the form $(\l x M ~ N)$] nor an $\eta$-redex [i.e. a $\l$-term
of the form $\l x (M ~ x)$ where $x$ does not appear in $M$].
\end{itemize} 

The following result is well known (B\H{o}hm Theorem):

\begin{itemize}
\item[] {\it If $U,V$ are two distinct closed $\b \eta$-normal $\l$-terms then there is a
closed $\l$-term $W$ such that $(W ~ U) \simeq\sb{\b} T$ and $(W ~ V) \simeq\sb{\b} F$.}
\end{itemize}  

\begin{itemize}
  \item Let us recall that a $\l$-term
$M$ either has a {\it head redex} [i.e. $M=\l x_1 ...\l x_n ((\l x U ~ V) ~ V_1 ... V_m)$, the
head redex being $(\l x U ~ V)$], or is in {\it head normal form} [i.e. $M=\l x_1 ...\l x_n (x
~ V_1 ... V_m)$]. \item The notation $U \p V$ means that $V$ is obtained from $U$ by some head
reductions and we denote by $h(U,V)$ the length of the head reduction between $U$ and $V$. 
\item A $\l$-term is said {\it solvable} iff its head reduction terminates.
\end{itemize} 

The following results are well known :

\begin{itemize}
\item[] {\it - If $M$ is $\b$-equivalent to a head normal form then $M$ is solvable.}
\item[]Ê{\it - If $U \p V$, then, for any substitution $\s$, $\s(U) \p \s(V)$,
and $h(\s(U),\s(V))$=h(U,V). \\
In particular, if for some  substitution $\s$, $\s(M)$ is
solvable, then $M$ is solvable.}
\end{itemize}  
 
\section{Numeral systems}

\begin{itemize}
\item A {\it numeral system} is an infinite sequence of different closed $\b \eta$-normal
$\l$-terms $\bf d$ \rm = $d_0 , d_1 ,..., d_n , ...$.
\item Let {\bf d} be a numeral system.
\begin{itemize}
\item A closed $\l$-term $S_d$ is called {\it Successor} for $\bf d$ iff :
 \begin{center}
$(S_d ~ d_n) \simeq\sb{\b} d_{n+1}$ for all $n \in \N$.
\end{center} 

\item A closed $\l$-term $P_d$ is called {\it Predecessor} for $\bf d$ iff :
 \begin{center}
$(P_d ~ d_{n+1}) \simeq\sb{\b} d_n$ for all $n \in \N$.
\end{center}

\item A closed $\l$-term $Z_d$ is called {\it Zero Test} for $\bf d$ iff :
 \begin{center}
$(Z_d ~ d_0) \simeq\sb{\b} T$ \\
and\\
$(Z_d  ~ d_{n+1}) \simeq\sb{\b} F$ for all $n \in \N$.
\end{center}
\end{itemize}

\item A numeral system is called {\it adequate} iff it possesses closed $\l$-terms for
Successor,  Predecessor, and Zero Test.
\end{itemize}

{\bf Examples of adequate numeral systems}\\

{\it 1) \so{The Barendregt numeral system}} \\
For each $n \in \N$, we define the Barendregt
integer $\sou{n}$ by : $\sou{0}=I$ and $\sou{n+1}=<F,\sou{n}>$. \\
It is easy to check that 
\begin{itemize}
\item[]Ê$\sou{S}=\l x < F , x>$, 
\item[] $\sou{P}=\l x(x ~ F)$, 
\item[] $\sou{Z}=\l x (x ~ T)$. 
\end{itemize}
are respectively $\l$-terms for Successor, Predecessor, and Zero Test. $\Box$ \\

{\it 2) \so{The Church numeral system}} \\
For each $n \in \N$, we define the Church integer
$\so{n}=\l f \l x (f ( f...(f ~ x)...))$ ($f$ occurs $n$ times).\\
It is easy to check that
\begin{itemize}
\item[] $\so{S}=\l n \l f \l x(f ~ (n ~ f ~ x))$, 
\item[] $\so{P}=\l n(n ~ U ~ <\so{0},\so{0}> ~ T)$ where $U=\l a<(\so{s}~ (a ~ T)),(a ~ F)>$, 
\item[] $\so{Z}=\l n (n ~ \l x F ~ T)$.  
\end{itemize}
are respectively $\l$-terms for Successor, Predecessor, and Zero Test. $\Box$ \\ 
 
Each numeral system can be naturally considered as a coding of integers into $\l$-calculus
and then we can represent total numeric functions as follows.
\begin{itemize}
\item A total numeric function $\phi : \N^p \f \N$ is {\it $\l$-definable} with respect to a 
numeral  system {\bf d} iff \begin{center}
$\e$ $F_{\phi}$ $\q$ $n_1,...,n_p \in \N$ $(F_{\phi} ~ d_{n_1}... d_{n_p}) \simeq\sb{\b}
d_{\phi(n_1,...,n_p)}$ 
\end{center}
\end{itemize}

The Zero Test can be considered as a function on integers. Indeed :\\

{\bf Lemma 1} {\it A numeral system {\bf d} has a $\l$-term for Zero Test iff the
function $\phi$ defined by : $\phi(0)=0$ and $\phi(n)=1$ for every $n \geq 1$ is
$\l$-definable with respect to {\bf d}.}  \\

{\bf Proof} It suffices to see that $d_0$ and $d_1$ are distinct $\b \eta$-normal $\l$-terms.
$\Box$ \\

H. Barendregt has shown in [1] that :\\

{\bf Theorem 1} {\it A numeral system {\bf d} \it is adequate iff all total recursive functions are
$\l$-definable with respect to} {\bf d}.   
 
\section{Some results on numeral systems}

{\bf Theorem 2} \footnote {This Theorem is the exercise 6.8.21 of Barendregt's book
(see [1]). We give here a proof based on the techniques developed by J.-L. Krivine in [3].}
{\it There is a numeral system with Successor and  Predecessor but without Zero Test.}  \\

{\bf Proof} For every $n \in \N$, let $a_n = \l x_1...\l x_n I$.\\
It is easy to check that the $\l$-terms $S_a = \l n \l x n$ and $P_a = \l n (n ~
I)$ are $\l$-terms for Successor and Predecessor for {\bf a}.\\
Let $\n,x,y$ be different variables. \\
If {\bf a} possesses a closed $\l$-term $Z_a$ for Zero Test, then : 
 \begin{center}
$(Z_a ~ a_{n} ~ x ~ y) \simeq\sb{\b} \cases { x &if $n=0$ \cr y &if $n \geq 1$
\cr}$ 
\end{center}
and
 \begin{center}
$(Z_a ~ a_{n} ~ x ~ y) \p \cases { x &if $n=0$ \cr y &if $n \geq 1$
\cr}$ 
\end{center}
Therefore  $(Z_a ~ \n ~ x ~ y)$ is solvable and its head normal form does not begin with
$\l$. \\
We have three cases to look at :  
\begin{itemize} 
\item   $(Z_a ~ \n ~ x ~ y) \p (x ~ u_1...u_k)$, then $(Z_a ~ a_1 ~ x ~ y) \not \p y$.
\item   $(Z_a ~ \n ~ x ~ y) \p (y ~ u_1...u_k)$, then $(Z_a ~ a_0 ~ x ~ y) \not \p x$. 
\item   $(Z_a ~ \n ~ x ~ y) \p (\n ~ u_1...u_k)$, then $(Z_a ~ a_{k+2} ~ x ~ y) \not \p y$.
\end{itemize}
Each case is impossible.  $\Box$ \\

{\bf Theorem 3} {\it There is a numeral system with Successor and Zero Test but without
Predecessor.}  \\

{\bf Proof} Let $b_0 = < T , I >$ and for every $n \geq 1$,  $b_n = <  F , a_{n-1} >$.\\ 
It is easy to check that the $\l$-terms $S_b = \l n < F , ((n ~ T) ~ a_0 ~ \l x (n~ F)) >$ and
$Z_b = \l n (n ~ T)$ are $\l$-terms for Successor and Zero Test for {\bf b}.\\ If {\bf b}
possesses a closed $\l$-term $P_b$ for Predecessor, then the $\l$-term $P'_b = \l n (P_b ~ < F
, n > ~ T)$ is a $\l$-term for Zero Test for {\bf a}. A contradiction. $\Box$ \\

{\bf Remarks}\\
1) Let $b'_0=b_1$, $b'_1=b_0$, and for every $n \geq 2$, $b'_n=b_n$. It is easy to check that
the numeral system {\bf b$'$} does not have $\l$-terms for Successor, Predecessor, and Test for
Zero.\\
2) The proofs of Theorems 1 and 2 rest on the fact that we are considering sequences of
$\l$-terms with a strictly increasing order (number of abstractions). Considering sequences
of $\l$-terms with a strictly increasing degree (number of arguments) does not work as well.
See the following example.\\  We define $\tilde{0}=I$ and for each $n \geq 1$, $\tilde{n}=\l x
(x ~ x...x)$ ($x$ occurs $n+1$ times). \\ Let  \begin{itemize} \item[]Ê$\tilde{S}=\l n \l x (n
~ x ~ x)$,  \item[] $\tilde{Z}=\l n (n ~ A ~ I ~ I ~ T)$ where $A=\l x\l y (y ~ x)$,  \item[]
$\tilde{P}=\l n \l x( n ~ U ~ F)$  where $U=\l y (y ~ V ~ I)$  and $V=\l a \l b \l c \l d (d ~
a ~ (c ~ x))$ \end{itemize}  It is easy to check that $\tilde{S}$, $\tilde{Z}$, and
$\tilde{P}$ are respectively $\l$-terms  for Successor, Zero Test, and Predecessor. 
$\Box$ \\ 

{\bf Definitions}
\begin{itemize}
\item We denote by $\L^{0}$ the set of closed $\l$-terms and by $\L^{1}$ the set of the
infinite sequences of closed normal $\l$-terms. It is easy to see that $\L^{0}$ is countable
but $\L^{1}$ is not countable.
\item For every finite sequence of $\l$-terms $U_1,U_2,...,U_n$
we denote by $<U_1,U_2,...,U_n>$ the $\l$-term $<...<<I,U_1>,U_2>,...,U_n>$. 
\item Let  {\bf U}=$U_1,U_2,...$ be a sequence of normal closed $\l$-terms. A closed 
 $\l$-term $A$ is called {\it generator} for {\bf U} iff :
\begin{center}
$(A ~ I) \simeq\sb{\b} U_1$ \\
and\\
$(A ~ <U_1,U_2,...,U_n>) \simeq\sb{\b} U_{n+1}$ for every $n \geq 1$ \\
\end{center}
\end{itemize}

{\bf Lemma 2} {\it There is a sequence of normal closed $\l$-terms without generator.}  \\ 

{\bf Proof } If not, let $\ph$ be a bijection between $\L^{0}$ and $\N$ and $\Ph$ the function
from
 $\L^1$ into $\L^0$ defined by: $\Ph({\bf U})$ is the generator $G_{{\bf U}}$ such that
$\ph(G_{{\bf U}})$ is minimum. It is easy to check that $\Ph$ is a one-to-one mapping.  A
Contradiction.  $\Box$ \\

{\bf Theorem 4} {\it There is a numeral system with Predecessor and Zero Test but
without Successor.}  \\

{\bf Proof} Let {\bf e} be a sequence of normal closed $\l$-terms without generator.\\ 
Let $c_0 = I$ and for every $n \geq 1$,  $c_n = <  c_{n-1} , e_n
>$.\\  It is easy to check that the $\l$-terms $P_c = \l n ( n ~ T)$ and $Z_c = \l n (n ~~
\l x\l y I ~~ T ~~ F ~~ T)$ are $\l$-terms for Predecessor and Zero Test for {\bf c}.\\
If {\bf c} possesses a closed $\l$-term $S_c$ for Successor, then the $\l$-term $S'_c = \l n
(S_c ~  n  ~ F)$ is a generator for {\bf e}.
A Contradiction.  $\Box$ \\

The result of H. Barendregt (Theorem 1) means that, for a numeral system, it suffices to
represent three particular functions in order to represent all total recurcive functions. We
have proved that these three particular functions are independent. We think it is, at least,
necessary to have three functions as is mentioned below : \\

{\bf Conjecture} {\it There are no total recursive functions $f,g : \N \f \N$ such
that : for all numeral systems {\bf d}, $f,g$ are $\l$-definable iff all total recursive
functions are $\l$-definable with respect to {\bf d}}.  \\

If we authorize the binary functions we obtain the following result :\\

{\bf Theorem 5} {\it There is a binary total function $k$ such that for all numeral systems
{\bf d}, $k$ is $\l$-definable iff all total recursive functions are $\l$-definable with
respect to {\bf d}}.\\

{\bf Proof} Let $k$ the total binary function defined by :
\begin{center}
$k(n,m) = \cases { n+1 &if $m=0$ \cr |n-m| &if $m \not = 0$ \cr}$
\end{center}

It suffices to see that :

\begin{center}
$k(n,n)= \cases {1 &if $n=0$ \cr 0 &if $n \not = 0$ \cr}$,
\end{center}

\begin{center}
$k(n,0)=n+1$,
\end{center}

\begin{center}
$k(n,1)=n-1$ if $n \not = 0$. $\Box$\\ [0,5cm]
\end{center}

{\bf Acknowledgement.} We wish to thank Mariangiola Dezani and Ren\'e David for helpful
discussions. We also thank Enrico Tronci and N\H{o}el Bernard for their help in the writing of
this paper.

\end{document}